\pdfoutput=1
\documentclass[11pt]{amsart}
\usepackage{graphicx,hyperref,amsmath,amsthm,bm,float}
\usepackage[numbers]{natbib}
\usepackage{tikz}
\usetikzlibrary{babel,cd}
\usepackage{microtype,todonotes}
\usepackage[english]{babel}
\usepackage[a4paper,text={16.5cm,25.2cm},centering]{geometry}
\newtheorem{thm}{Theorem}[section]
\newtheorem{coro}[thm]{Corollary}
\newtheorem{prop}[thm]{Proposition}
\newtheorem{cons}[thm]{Construction}
\newtheorem{lemma}[thm]{Lemma}
\newtheorem{ques}[thm]{Question}
\theoremstyle{definition}
\newtheorem{defn}[thm]{Definition}
\theoremstyle{remark}
\newtheorem{rem}[thm]{Remark}

\title{Uniform Foliations with Reeb components}
\author{Joaqu\'in Lema}
\address{Centro de Matem\'atica, UdelaR, Uruguay}
\email{joalema@cmat.edu.uy}
\thanks{The author was supported by CAP grant for Master's students at UdelaR}

\begin{document}

\begin{abstract}
A foliation on a compact manifold is uniform if each pair of leaves of the induced foliation on the universal cover are at finite Hausdorff distance from each other. We study uniform foliations with Reeb components. We give examples of such foliations on a family of closed $3-$manifolds with infinite fundamental group. Furthermore, we prove some results concerning the behavior of a uniform foliation with Reeb components on general $3-$manifolds.
\end{abstract}

\maketitle

\section{Introduction}

Consider a foliation $\mathcal{F}$ on a compact Riemannian $3-$manifold $M$. This foliation lifts to a foliation $\widetilde{\mathcal{F}}$ on the universal cover $\widetilde{M}$. We will say that $\mathcal{F}$ is \emph{uniform} if any pair of leaves of $\widetilde{\mathcal{F}}$ are at finite Hausdorff distance from each other. 

A lot can be said about a uniform foliation if we further assume $\mathcal{F}$ to be Reebless (see for example \cite{FP20}, \cite{Th97}). In this paper, we will focus on the opposite case. More precisely, we will study the following question posed by S. Fenley and R. Potrie on the article \cite[Question 1]{FP20}:

\begin{ques}
\label{ques:FP}
If $\mathcal{F}$ is uniform in $M$ with infinite fundamental group, does it follow that $\mathcal{F}$ is also Reebless?
\end{ques}

Our first result will be to give a negative answer to the question (see Section \ref{section:constructions}):

\begin{thm}
\label{thm:constructions}
For every $l,m \in \mathbb{N}$ and every choice $M_1,\ldots,M_l$ of $3-$manifolds with finite fundamental group, there exists a uniform foliation with Reeb components on $M = \left( \#_{i=1}^m S^1 \times S^2 \right) \# \left( \#_{i=1}^l M_i \right)$.
\end{thm}

One can arrange the foliations given by Theorem \ref{thm:constructions} to be $C^\infty-$smooth (see Remark \ref{rem:smoothness}). We will also give an example of uniform foliation with Reeb components on the solid torus, which is trivial on the boundary. However, these $3-$manifolds are ``small'' in the sense that they are the only ones not admitting an immersed essential closed surface of genus $g \geq 1$ (see Proposition \ref{prop:manifoldsadmitingessential}). In this paper, we will say that an immersed surface is \emph{essential} if the immersion induces an injective morphism from the fundamental group of the surface to $\pi_1 (M)$.

We can say the following for the rest of $3-$manifolds:

\begin{thm}
\label{thm:main}
Let $\mathcal{F}$ be a uniform foliation on a compact $3-$manifold $M$ admitting an immersion $i: \Sigma \rightarrow M$, from a closed surface $\Sigma$ of genus $g \geq 1$, such that $i_* : \pi_1 (\Sigma) \rightarrow \pi_1 (M)$ is injective. Then $i(\Sigma)$ must intersect the set of Reeb components.
\end{thm}

The idea of the proof is to put the immersion in general position with respect to the foliation, and then apply Poincar\'e-Bendixson theory on the induced foliation on $\widetilde{\Sigma} \cong \mathbb{R}^2$. 

Theorem \ref{thm:main} and some remarks in Section \ref{sec:remarks} motivate us to modify Question \ref{ques:FP}: does an \emph{irreducible} $3-$manifold with infinite fundamental group admit a uniform foliation with Reeb components? 

\subsection*{Acknowledgments} 

I am profoundly grateful to my advisor Rafael Potrie for his patience and suggestions during every step of my thesis. This work could not be done without him. I would also like to thank Sergio Fenley for his comments on earlier versions of this paper that motivated the statement of Theorem \ref{thm:main}.

\section{Preliminaries}

\subsection{Foliations on $3-$manifolds:} We will be working on a compact $3-$manifold $M$ endowed with a $C^{\infty,0+}$ codimension one foliation $\mathcal{F}$. From now on, by foliation on a $3-$manifold we refer to a foliation of codimension one. We will assume some familiarity with foliation theory; the reader may find a comprehensive treatment of the subject in \cite{CC00,HH83}.

A Reeb component is a foliation of the solid torus so that the boundary is a leaf, and there is a circle worth of planar leaves in the interior spiraling towards the boundary torus. We will also call some quotient of this foliation a Reeb component. These are crucial on the study of foliations on $3-$manifolds by Novikov's celebrated theorem \cite{Nov65}. It says (among other things) that if $\mathcal{F}$ is orientable and transversely orientable, then it is Reebless if and only if every leaf $L$ is essential (as defined in the introduction). In this sense, for Reebless foliations the topology of a leaf is tied to the topology of the manifold.

Considering this fact is natural to ask if there exists some property of a foliation $\mathcal{F}$ with Reeb components ``reading'' the topology of the $3-$manifold. This does not seem plausible as every $3-$manifold admits a foliation with Reeb components. Furthermore, Thurston shows in \cite{Thu76} that every plane field on $M$ is homotopic to the tangent space of a foliation. His construction is local in nature, so we do not care about the global topology of $M$. The Reeb components play a key role because holes can easily be filled using them. See \cite[Section 8.5]{CC00} for a detailed treatment of this construction. 

However, we believe that if a compact $3-$manifold $M$ is ``big enough'', then it does not admit a \emph{uniform} foliation with Reeb components (see Section \ref{sec:remarks}). Particularly, Lemma \ref{lemma:FP} would tell us that every foliation on these manifolds must have leaves lifting to unbounded sets on the universal cover. This would be interesting because not much is known about the behavior of a foliation with Reeb components on a general $3-$manifold.

Suppose that $M$ is equipped with a Riemannian metric, and let us denote the universal cover of $M$ by $p: \widetilde{M} \rightarrow M$, equip $\widetilde{M}$ with the pullback metric.

\begin{defn}
Let $\mathcal{F}$ be a foliation on a compact manifold $M$ and denote the lifted foliation on $\widetilde{M}$ by $\widetilde{\mathcal{F}}$. We will say that $\mathcal{F}$ is \emph{uniform} if any pair of leaves of $\widetilde{\mathcal{F}}$ are at finite Hausdorff distance from each other.
\end{defn}

This notion was introduced by Thurston in \cite{Th97}, who further required $\mathcal{F}$ to be Reebless (see also \cite[Section 9.3]{Ca07}). We stick to the definition given by S. Fenley and R. Potrie in \cite{FP20}. 

\subsection{Turbulization} 

To prove Theorem \ref{thm:constructions} we will rely on a method for modifying foliations along a transverse curve known as \emph{turbulization}. 

Suppose that a $3-$manifold $M$ is endowed with an oriented and cooriented foliation $\mathcal{F}$ admitting an embedded closed curve $\tau$ transverse to the foliation. Let $\overline{N(\tau)}$ be the closure of an embedded tubular neighborhood of $\tau$, then the orientability and coorientability of $\mathcal{F}$ implies that if $N(\tau)$ is small enough, $\mathcal{F}|_{\overline{N(\tau)}}$ is homeomorphic to the product foliation by disks on the solid torus $\overline{N(\tau)}$.

Fix an identification of $\overline{N(\tau)}$ with $D^2 \times S^1$ sending the leaves of $\mathcal{F}|_{\overline{N(\tau)}}$ to $D^2 \times \{\cdot\}$. We will denote points in $D^2 \times S^1$ in cylindrical coordinates, by this we mean that $(r,\theta,z)$ represents a point whose projection to $D^2$ has polar coordinates $(r,\theta)$ and projecting to $z \in S^1$. Notice that the plane field tangent to the foliation by disks is given by the kernel of the $1-$form $\alpha_0 = dz$.

Now we will construct a new foliation $\mathcal{T}$ on $D^2 \times S^1$ with a Reeb component on the interior and coinciding with the foliation by disks on a neighborhood of the boundary torus. This will give us a new foliation $\mathcal{F}'$ on $M$ defined as $\mathcal{F}$ in the complement of $\overline{N(\tau)}$ and as $\mathcal{T}$ on $\overline{N(\tau)}$.

\begin{figure}[h]
\begin{center}
\includegraphics[scale=0.8]{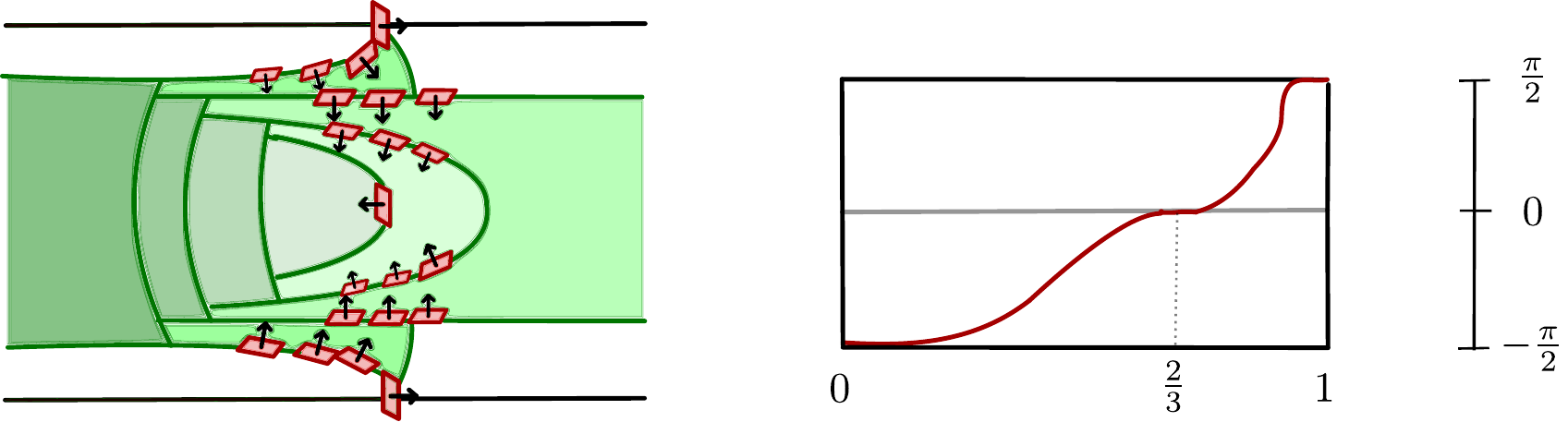}
\caption{At the left of the picture: some leaves of the foliation $\mathcal{T}$. At the right: the function $\lambda$ controlling the vector field normal to $\ker \omega$.}
\label{fig:turbulization}
\end{center}
\end{figure}

Take $\lambda : [0,1] \rightarrow [-\frac{\pi}{2},\frac{\pi}{2}]$ a smooth function such that it is a strictly increasing bijection restricted to the interval $[0,\frac{3}{4}]$, $\lambda (\frac{2}{3}) = 0$ and $\lambda|_{[\frac{3}{4}, 1]} = \frac{\pi}{2}$ (see the right of Figure \ref{fig:turbulization}). Now define the $1-$form:
$$\omega_{(r,\theta,z)} = \cos (\lambda (r) ) dr + \sin (\lambda(r)) dz.$$
This is well defined on $D^2 \times S^1$ if $\lambda^{(k)} (0) = 0$ for every $k \geq 1$. Using Frobenius' theorem, one can check that $\pi = \ker \omega$ is an integrable plane field. The foliation tangent to $\pi$ is the desired foliation $\mathcal{T}$ depicted on the left of Figure \ref{fig:turbulization}.

If we are careful with the choice of $\lambda$, it can be shown that the resulting foliation $\mathcal{F}'$ on $M$ is as regular as $\mathcal{F}$. The reader can find a detailed description of this construction in \cite[Example 3.3.11, vol. 1]{CC00}.

\subsection{General Position} 

Let $\mathcal{F}$ be a foliation on a $3-$manifold $M$ and $i: \Sigma \rightarrow M$ an immersion. We can always take a small perturbation of $i$ in order to suppose that $i(\Sigma)$ is ``as transverse as possible'' to the foliation. 

\begin{defn}
\label{definition:generalposition}
Let $\Sigma$ be a closed surface and $i: \Sigma \rightarrow M$ an immersion into a $3-$manifold $M$ endowed with a foliation $\mathcal{F}$. We will say that \emph{the immersion is in general position with respect to $\mathcal{F}$} if the following happens:
\begin{enumerate}
\item Except at a finite set of points $\{p_1,\ldots,p_l\} \subset \Sigma$, $i$ is transverse to $\mathcal{F}$.
\item The points $\{i(p_1),\ldots,i(p_l) \}$ lie on distinct leaves of $\mathcal{F}$.
\item For every $p_k$ and every submersion $\phi$ defined on a neighborhood around $i(p_k)$ locally defining the foliation \footnote{By this we mean a map $\phi : U \rightarrow \mathbb{R}$ such that the preimages of the regular values are disks contained on a leaf.} sending $i(p_k)$ to zero, there exist a neighborhood $U_k$ of $p_k$ such that $\phi \circ i \vert_{U_k}$ is topologically conjugated to $f_c$ or $f_s$, where $f_c (x,y) = x^2 + y^2$ and $f_s (x,y) = x^2 - y^2$. By this we mean that there exists a homeomorphism between $U_k$ and some neighborhood $U_k'$ of $0 \in \mathbb{R}^2$ such that the following diagram commutes:
\begin{center}
\begin{tikzcd}
U_k \arrow[d,"h"] \arrow[r,"(\phi \circ i)\vert_{U_k}"] &  \mathbb{R} \\
U_k' \arrow[ur,"f_\sigma"],
\end{tikzcd}
\end{center}
where $\sigma$ may be $c$ or $s$.
\end{enumerate}
\end{defn}

\begin{figure}[h!]
\begin{center}
\includegraphics[scale=0.7]{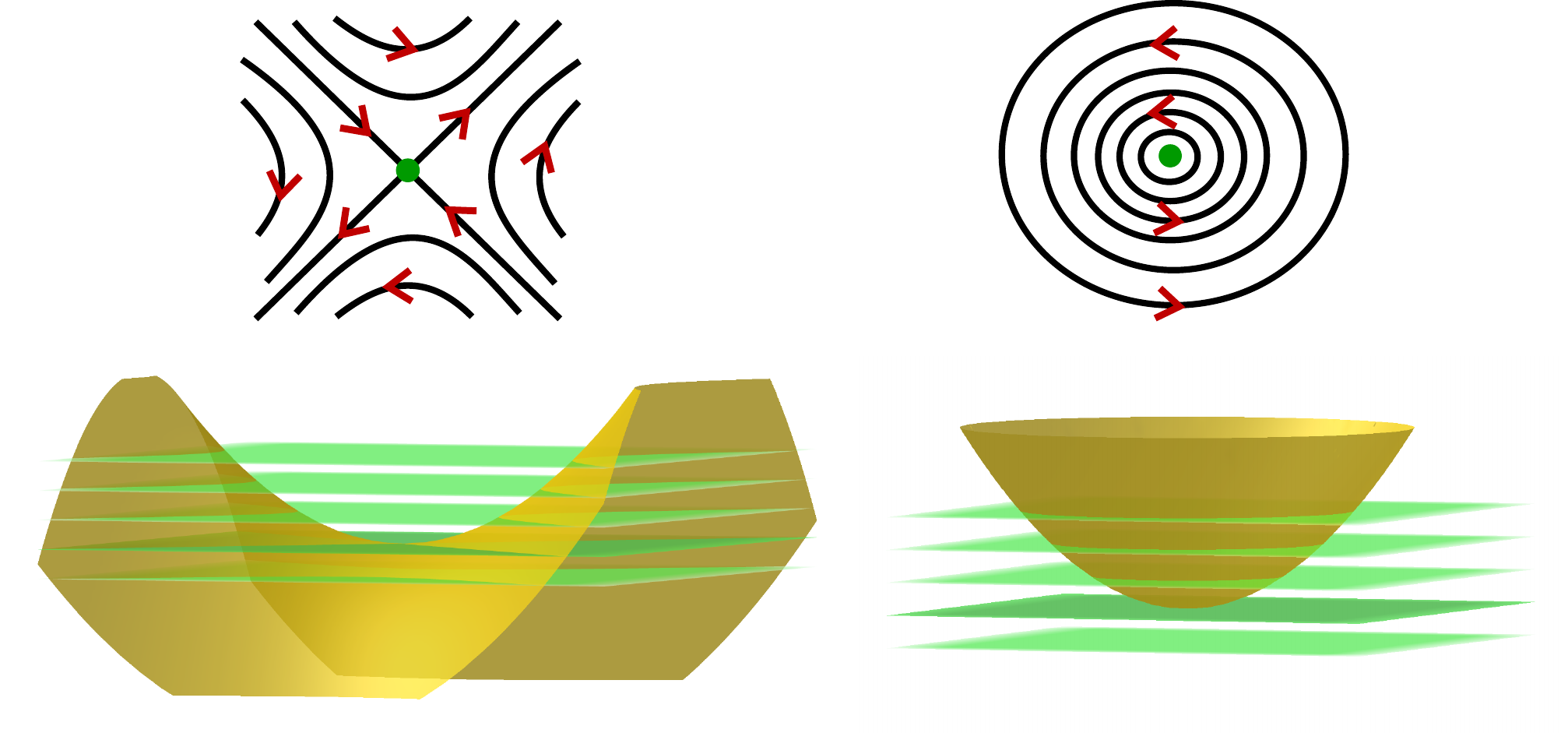}
\caption{At the bottom of the picture: tangencies of a surface in general position (in yellow) and the foliation (in green). At the top: the induced singular foliation on the surface. The orientation and coorientation of $\mathcal{F}$ induces an orientation on the singular foliation.}
\label{fig:genpos}
\end{center}
\end{figure}

Geometrically the last condition says that the tangencies between $\Sigma$ and the foliation look like a critical point of a Morse function on $\Sigma$ (up to homeomorphism), see figure 
\ref{fig:genpos}.

The following theorem dates back to Haefliger when the foliation $\mathcal{F}$ is sufficiently regular. It was generalized to the case of $C^0$ foliations by Solodov in \cite{Sol84}:

\begin{thm}
\label{thm:generalposition}
Let $i: \Sigma \rightarrow M$ be an immersion of a closed surface $\Sigma$ on a $3-$manifold $M$ endowed with a oriented and cooriented foliation $\mathcal{F}$. Then for every $\varepsilon > 0$ there exists an immersion $j : \Sigma \rightarrow M$ in general position with respect to $\mathcal{F}$ and $\varepsilon-$close to $i$ in the $C^0-$topology.
\end{thm}

Suppose that $\mathcal{F}$ is oriented and cooriented and $i: \Sigma \rightarrow M$ is an immersion in general position with respect to $\mathcal{F}$. The first condition of definition \ref{definition:generalposition} tells us that $\mathcal{F}$ induces a foliation on $\Sigma \setminus \{p_1,\ldots, p_l\}$ whose leaves are intersections of leaves of $\mathcal{F}$ and $i(\Sigma)$. The third condition gives us a model neighborhood around a singularity; more precisely they look like saddles in the case of $f_s$ or centers in the case of $f_c$ (see figure \ref{fig:genpos}). This local models gives us what is called a \emph{singular foliation with Morse singularities} on $\Sigma$. The orientation of $\mathcal{F}$ induces an orientation on this singular foliation (this is an orientation outside its singular points).

The second condition tells us that a separatrix cannot join two distinct singularities (see below for a definition of separatrix).

\subsection{Singular Foliations of Morse-type}

In this section, we will fix some notations that we are going to use throughout this text. Let $\Sigma$ be a surface endowed with an oriented singular foliation with Morse singularities $\mathcal{G}$. Suppose that $p$ is some non-singular point, we will denote the leaf passing through $p$ by $L_p$. If $L_p$ is non-compact (ie. it is not a circle) then $p$ separates $L_p$ into two components $L_p^+$ and $L_p^-$, where $L_p^+$ is composed of points greater than $p$ with respect to the order imposed by the orientation and $L_p^-$ of points smaller than $p$.

If $L$ is some non-compact leaf, we will say that the $\omega-$limit of the leaf is the set $\omega (L) = \bigcap_{p \in L} \overline{L_p^+}$. Analogously we define the $\alpha-$limit of the leaf as $\alpha (L) = \bigcap_{p\in L} \overline{L_p^-}$. A leaf $L$ is a \emph{separatrix} if the $\alpha$ or $\omega-$limit of $L$ is a singularity. 

Let $C$ be a union of singularities and separatrixes $S_i$ such that the $\alpha$ and $\omega-$limit of $S_i$ is a singularity in $C$. This set defines a directed graph with one vertex for each singularity in $C$ and one edge for each $S_i$, with the orientation induced by the foliation. We will say that $C$ is a \emph{closed graph} if there exists a closed path on the graph which travels through every edge only once (respecting the orientation).

We will end up studying Morse-type singular foliations on $\mathbb{R}^2$. In this context, we can apply the classical Poincar\'e-Bendixson theorem (see for instance \cite[Theorem 1.8]{PdM12}):

\begin{thm}[Poincar\'e-Bendixson]
Let $\mathcal{G}$ be a Morse type singular foliation on $\mathbb{R}^2$ and $L$ a leaf with compact closure. Then $\alpha (L)$ and $\omega (L)$ can be a saddle singularity, a closed leaf, or a closed graph.
\end{thm}

\subsection{Essential immersions}

Let $M$ be some closed $3-$manifold. We will say that a $2-$sided embedded closed surface $\Sigma$ is \emph{compressible} if there exists some embedded disk $D$ on $M$ such that $D \cap \Sigma = \partial D$, and $\partial D$ is not homotopically trivial in $\Sigma$. We will call $D$ a \emph{compressing disk} for $\Sigma$. A $2-$sided embedded surface $\Sigma$ of genus $g \geq 1$ is \emph{incompressible} if there are no compressing disks.

If a $2-$sided embedded surface $\Sigma$ is compressible, we can do surgery on a compressing disk $D$ in order to obtain a simpler surface \cite[Section 1.2]{Hat}. This operation preserves the homology class $[\Sigma]\in H_2 (M,\mathbb{Q})$. Doing finitely many surgeries on $\Sigma$, we obtain an embedded $2-$sided surface, whose connected components are incompressible surfaces or spheres. The sum of the homology class of the connected components is the homology class of $\Sigma$.

The disk theorem \cite[Theorem 3.1]{Hat} tells us that a $2-$sided embedded surface $\Sigma$ is incompressible if and only if the embedding $i: \Sigma \rightarrow M$ induces an injective morphism $i_* : \pi_1 (\Sigma) \rightarrow \pi_1 (M)$, i.e. the surface is essentially embedded.

Thanks to the virtual Haken conjecture proved by Agol in \cite{Ago13}, we can say exactly which closed $3-$manifolds admit a essentially embedded closed surfaces of genus $g\geq 1$:

\begin{prop}
\label{prop:manifoldsadmitingessential}
A closed $3-$manifold $M$ admits an essentially immersed surface $\Sigma$ of genus $g \geq 1$ if and only if some factor of the prime decomposition of $M$ is irreducible with an infinite fundamental group.
\begin{proof}
The ``if'' part is a direct consequence of the virtual Haken conjecture. Suppose $P$ is a factor of the prime decomposition, which is irreducible with an infinite fundamental group. Then it has a finite cover $\widehat{P}$ which admits an embedded incompressible surface of genus $g \geq 1$. Projecting this surface to $P$ and avoiding the balls of the connected sums, we obtain an essential immersion on $M$.

To see the ``only if'' part, it is enough to show that the fundamental group of the remaining $3-$manifolds does not admit a subgroup isomorphic to the fundamental group of a closed surface. No factor of the prime decomposition of $M$ is irreducible with an infinite fundamental group. Therefore, factors can be irreducible with a finite fundamental group or $S^1 \times S^2$ (prime but not irreducible). The fundamental group of such $3-$manifold must be a free product of infinite cyclic and finite groups.

If the fundamental group of a closed surface $\Gamma$ were isomorphic to a subgroup of such a group, then Kurosh subgroup Theorem \cite[Section 1.5.5]{Ser80} tells us that $\Gamma$ is isomorphic to a free product of infinite cyclic groups and finite groups. However, the fundamental group of a closed surface is freely indecomposable; this means that it is not isomorphic to a free product of two nontrivial groups (this follows from Stallings end theorem \cite[Section 4.A.6]{Sta72}). Therefore it would be isomorphic to a cyclic group or a finite group. Both options are impossible.
\end{proof}
\end{prop}

We remark that the last proof shows that the only closed $3-$manifolds which do not admit an essential immersion are those where we can use Theorem \ref{thm:constructions}.

\section{Some remarks}
\label{sec:remarks}

We begin with some general facts about uniform foliations \emph{with Reeb components} on $3-$manifolds. We will assume that $M$ is equipped with some Riemannian metric; this induces a metric on the universal cover $p : \widetilde{M} \rightarrow M$.

The following was observed by S. Fenley and R. Potrie in \cite{FP20}, we include the proof for completeness:

\begin{lemma}
\label{lemma:FP}
Let $\mathcal{F}$ be a uniform foliation with Reeb components on a compact manifold $M$, then every lift of a leaf to the universal cover $\widetilde{M}$ has compact closure. In particular the inclusion of a leaf $i: L \rightarrow M$ induces a morphism $i_* : \pi_1 (L) \rightarrow \pi_1 (M)$ with finite image.
\begin{proof}
It is enough to show that a boundary leaf of some Reeb component lifts to a compact leaf because any pair of leaves of the foliation $\widetilde{\mathcal{F}}$ are at finite Hausdorff distance from each other. So assume by contradiction that the lift of a boundary leaf of some Reeb component $R$ is noncompact. Let $q \in R$ and $\gamma \subset R$ be a loop at $q$ homotopic to the core of the Reeb component, the homotopy class of this curve is nontrivial; otherwise, $R$ would lift to a compact set. Choose some lift of $\gamma$ starting at $\widetilde{q} \in p^{-1} (q)$, we claim that $d_{\widetilde{M}} (\widetilde{q}, \gamma^n (\widetilde{q}) ) \to_n \infty$.

One way to see this is to use that the map $\alpha \in \pi_1 (M) \to \alpha (\widetilde{p})$ is a quasi-isometry (by Milnor-Svarc lemma). This fact implies that the distance between a lift of a plane on the interior of the Reeb component and $\gamma^n (\widetilde{q})$ goes to infinity with $n$, contradicting the uniform condition.

To see that the inclusion of a leaf $i: L \rightarrow M$ induces a morphism with finite image, suppose by contradiction that $\# (i_* \pi_1 (L)) = \infty$. Then, the cardinality of the stabilizer of a connected component $\widetilde{L}$ of $p^{-1}(L)$ is infinite. In particular, the orbit of every $p \in \widetilde{L}$ is infinite. But the action is proper, so $\widetilde{L}$ escapes every compact set, contradicting the last paragraph.
\end{proof}
\end{lemma}

This lemma motivates us to look for counterexamples in foliations such that every leaf of $\widetilde{\mathcal{F}}$ has compact closure. The following gives us an criterion for a foliation to verify this condition:

\begin{lemma}
\label{lemma:construction}
Let $\mathcal{F}$ be a foliation on a compact manifold $M$. Suppose that there exists $A_1,\ldots,A_k$ compact leaves such that every connected component of the complement of $\bigcup_{i=1}^k p^{-1} (A_i)$ has compact closure, then $\mathcal{F}$ is uniform.
\begin{proof}
Take a leaf $L$ of $\widetilde{\mathcal{F}}$ different from a lift of some $A_i$.
Then some connected component of the complement of $\bigcup_{i=1}^k p^{-1} (A_i)$ contains $L$ because the boundaries of these regions are composed of leaves. Therefore, the leaf has compact closure because each one of these regions is compact.
\end{proof}
\end{lemma}

This easy remark motivated our first counterexamples to Question \ref{ques:FP} (see Construction \ref{cons:S1S2} and \ref{cons:borromean}). However, the reader may notice that the existence of these compact leaves imposes conditions in the topology of the $3-$manifold $M$. Lemma \ref{lemma:FP} allows us to use compression disks on compact leaves until we get spheres on $M$. These spheres must bound some topology in order for us to use Lemma \ref{lemma:construction}.

The following is a consequence of this idea:

\begin{lemma}
\label{lemma:tori}
Let $\mathcal{F}$ be an oriented and cooriented uniform foliation with Reeb components on a closed, irreducible $3-$manifold $M$ with an infinite fundamental group. Then every compact leaf is a torus, bounding a solid torus in one of its sides.
\begin{proof}
The inclusion $i$ of a compact leaf $L$ induced the zero morphism $i_* : \pi_1 (L) \rightarrow \pi_1 (M)$, because the fundamental group of a closed, irreducible $3-$manifold is torsion free. We can compress the surface $L$ using the disk theorem until we get a union of spheres $S$. The compression operation preserves the class $[L] \in H_2 (M;\mathbb{Q})$. However $S$ is composed of  spheres and $M$ is irreducible, therefore $[L] = [S] = 0$. 

This implies that the leaf $L$ cannot admit a closed transversal $\tau$, because if it did then the intersection product between $[\tau] \in H_1 (M;\mathbb{Q})$ and $[L] \in H_2 (M;\mathbb{Q})$ would be non trivial and $[L] \neq 0$. A theorem of Goodman (see \cite[Theorem 6.3.5, vol. 1]{CC00}) implies that $L$ must be a torus. To see that it bounds a solid torus, do surgery on $L$ with some compressing disk to get a sphere. This sphere bounds a ball on one side, so we obtain a solid torus bounded by $L$ by reversing the compression operation. 
\end{proof}
\end{lemma}

\begin{figure}[h]
\begin{center}
\includegraphics[scale=1.2]{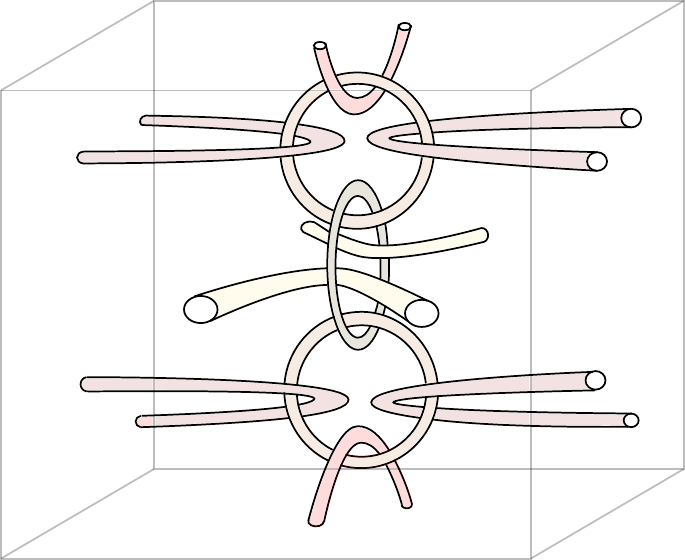}
\caption{Can we have a uniform foliation on $\mathbb{T}^3$ so that the Reeb components are the one depicted in this image?}
\label{fig:nightmaretori}
\end{center}
\end{figure}

Lemmas \ref{lemma:FP} and \ref{lemma:tori} are painting a strange picture. For instance, if the universal cover is $\mathbb{R}^3$ then the foliation is filling all the space, but every compact leaf cannot bound topology. An exceptional minimal set is not much worse because they must have compact closure. On the other hand, our Main Theorem \ref{thm:main} tells us that the set of Reeb components must intersect every essentially immersed surface. Figure \ref{fig:nightmaretori} shows a possible picture for the set of Reeb components on a uniform foliation on $\mathbb{T}^3$.

\begin{ques}
\label{ques}
If $\mathcal{F}$ is uniform in an irreducible $3-$manifold $M$ with an infinite fundamental group, does it follow that $\mathcal{F}$ is also Reebless?
\end{ques}

\section{Constructions}
\label{section:constructions}

We begin with an explicit counterexample to Question \ref{ques:FP} on $S^1 \times S^2$. To do so, consider the characterization of $S^1 \times S^2$ given by the following remark:

\begin{rem}
\label{rem:charactS1S2}
The manifold $S^1 \times S^2$ is homeomorphic to $M$, which is constructed as follows. Start with the solid torus $S^1 \times D^2$ and a circle $K$, which is the boundary of an embedded disk $D$. Identify a tubular neighborhood of the knot $K$ with $S^1 \times D^2$ in such a way that the longitude $S^1 \times \{\cdot\}$ bounds an embedded disk in the complement of the neighborhood. Drilling this tubular neighborhood, we obtain a compact $3-$manifold $N$ with two boundary components homeomorphic to $S^1 \times S^1$. Identifying both boundary components via the identity map, we obtain the closed $3-$manifold $M$. 

To prove this, notice that $S^1 \times S^2$ is homeomorphic to the quotient of $\mathbb{R}^3 \setminus \{0\}$ induced by the properly discontinuous $\mathbb{Z}-$action generated by the map $f(x) = \frac{x}{2}$. If $S^2$ is the standard sphere, the set bounded by $S^2$ and $f(S^2)$ is a fundamental domain. Now consider a torus $S$ resulting from $S^2$ by adding a small handle. Notice that $f(S)$ is contained in the solid torus bounded by $S$ and containing $0$. Furthermore, the set $N$ bounded by $S$ and $f(S)$ is also a fundamental domain for the action. Therefore, $S^1 \times S^2$ is homeomorphic to $N$, where we identify both boundary components via the map $f$. This identification coincides with the identification defining $M$.
\end{rem}

We will be referring to a curve $K$ on a $3-$manifold which is the boundary of an embedded disk as an \emph{unknot}. We will always equip an unknot with a frame such that the longitude of the tubular neighborhood is the boundary of some embedded disk in the complement of the neighborhood. 

\begin{figure}[t]
\begin{center}
\includegraphics[scale=0.8]{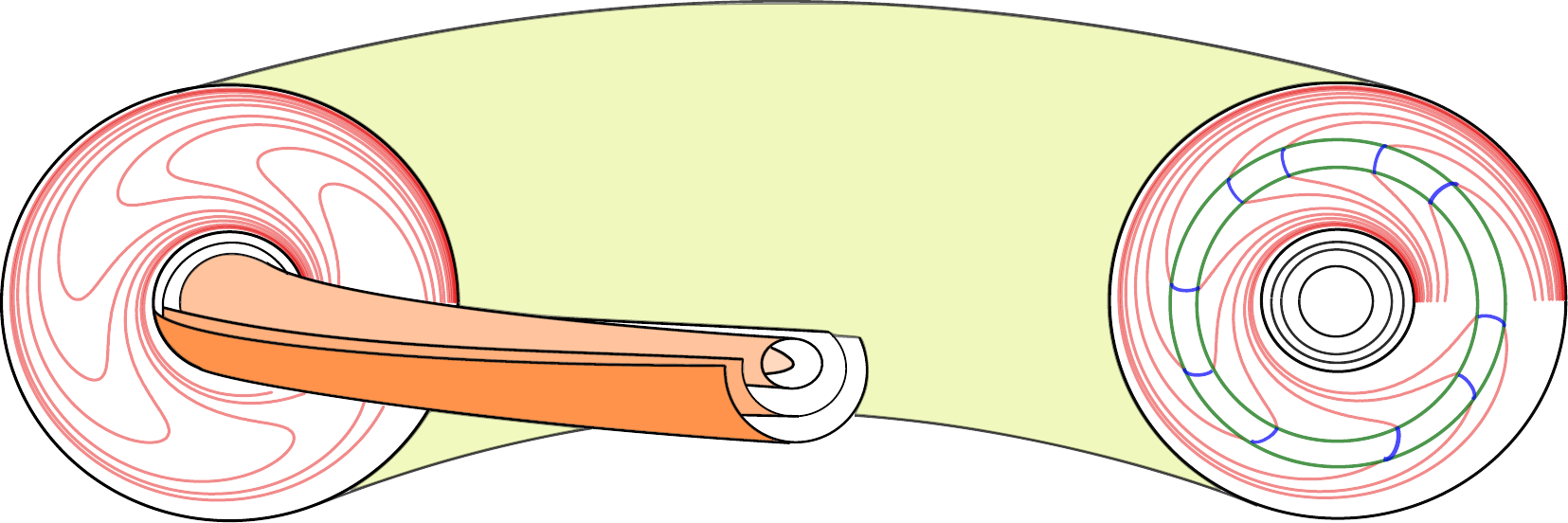}
\caption{The foliation on $N$ before turbulizing along $\tau$.}
\label{fig:folS1S2}
\end{center}
\end{figure}

The following construction is motivated by Lemma \ref{lemma:construction}:

\begin{cons}
\label{cons:S1S2}
We will construct a uniform foliation with Reeb components $\mathcal{F}$ on $S^1 \times S^2$ which admits a transversal unknot.

Thanks to Remark \ref{rem:charactS1S2}, it is enough to find a foliation tangential to the boundary of the manifold $N$ obtained from $S^1 \times D^2$ by drilling a tubular neighborhood of an unknot. To see this, notice that $S^1 \times S^2$ is obtained from $N$ by identifying the boundary components. The foliation on $N$ projects to a foliation of $S^1 \times S^2$ with a toric leaf $S$ coming from the boundary components of $N$. We saw on Remark \ref{rem:charactS1S2} that the lifts of $S$ to the universal cover $\mathbb{R}^3 \setminus \{0\}$ bounds a compact fundamental domain. Therefore this foliation falls under the hypothesis of Lemma \ref{lemma:construction}. In particular, it is uniform.

We may construct a foliation tangential to the boundary on $N$ as follows: let $D^2$ be the unit disk on $\mathbb{R}^2$ and $A \subset D^2$ be the annuli with inner and outer radii being $\frac{1}{2}$ and $1$ respectively. Equip $A$ with a Reeb component and consider the product foliation induced on $S^1 \times A \subset S^1 \times D^2$. We can complete this foliation to $S^1 \times D^2$ by adding a Reeb component on $S^1 \times \left( D^2 \setminus A \right)$. Notice that the core of a Reeb annulus is transversal to the foliation; this induces a transversal $\tau$ to the foliation on the solid torus bounding a disk. Turbulizing along $\tau$ and drilling the newly generated Reeb component, we obtain the desired foliation on $N$ (see Figure \ref{fig:folS1S2}).
\end{cons}

\begin{rem}
\label{rem:smoothness}
The constructed foliation can be made $C^\infty$. We just have to choose the Reeb annulus and the turbulization process in such a way that the resulting foliation on $N$ has $C^\infty-$trivial holonomy on the boundary (see \cite[Definition 3.4.1]{CC00}). The reader may check that each construction on this section can be made $C^\infty$. We will not pay attention to this.
\end{rem}

Now we will describe a procedure that will allow us to start with a pair of foliated $3-$manifolds $M,N$ and construct an explicit foliation on its connected sum. Furthermore, if we look at $M$ and $N$ with balls removed inside $M\# N$, the foliation on the connected sum coincides with the original pair of foliations outside the boundary of the drilled balls. The idea comes from the following remark:

\begin{rem}
\label{rem:connectedsums}
Let $M_1,M_2$ be a pair of $3-$manifolds. Consider an unknot $K_i$ on each $M_i$ ($i=1,2$) contained in a certain $3-$ball $B_i \subset M_i$ (for example, a tubular neighborhood of the spanned disk). Drilling a tubular neighborhood of $K_i$ contained on $B_i$ from $M_i$, we obtain $N_i$. Identifying this tubular neighborhood with $S^1 \times D^2$ via the standard framing, we get an identification between the boundary component of $N_i$ coming from $K_i$ and $S^1 \times S^1$.

We claim that if we start from $N_1 \cup N_2$ and we identify the boundaries of the drilled tubular neighborhoods of $K_1$ and $K_2$ via the map $f: S^1 \times S^1 \rightarrow S^1 \times S^1$ defined as $f(x,y) = (y,x)$, we obtain a $3-$manifold $M$ homeomorphic to the connected sum $M_1 \# M_2$. To prove this claim, notice that if we start with the balls $B_i$, we drill an unknot $K_i$ from each of them, and we sew the resulting sets along the boundary tori according to the map $f$, then we obtain $S^3$ with a pair of $3-$balls removed, i.e. $S^2 \times [0,1]$. Therefore $M \cong M_1 \# S^3 \# M_2$ as desired.
\end{rem}

\begin{lemma}
\label{lemma:connectedsum}
Let $M_1,M_2$ be compact $3-$manifolds endowed with foliations $\mathcal{F}_1,\mathcal{F}_2$. Suppose that there exists an unknot $K_i$ in $M_i$, transverse to $\mathcal{F}_i$ and contained in a ball $B_i$. Then there exist a foliation $\mathcal{F}$ on $M_1 \# M_2$ admitting a transversal unknot and coinciding with $\mathcal{F}_i$ on $M_i \setminus B_i$. Furthermore if $\mathcal{F}_1$ and $\mathcal{F}_2$ are uniform, then $\mathcal{F}$ is also uniform.
\begin{proof}
We start by turbulizing $\mathcal{F}_i$ along the transversal $K_i$. This process can be done by letting the foliation $\mathcal{F}_i$ fixed outside the ball $B_i$. Drilling the newly generated Reeb component, we obtain a $3-$manifold $N_i$ equipped with a foliation tangential to the component $S_i$ of $\partial N_i$ coming from the Reeb component. Identifying $S_i \subset \partial M_i$ as in Remark \ref{rem:connectedsums}, we obtain the desired foliation $\mathcal{F}$ on $M_1 \# M_2$. Notice that the leaf obtained by identifying $S_i$ is a torus $S$ such that the inclusion $i: S \rightarrow M_1 \# M_2$ induces the zero morphism on the fundamental group. Near this leaf, there are a lot of transversal unknots resulting from the turbulization process.

Now we will see that this construction preserves the uniform condition. Notice that as $\mathcal{F}_i$ admits a homotopically trivial transversal, it must have Reeb components by Novikov's Theorem (see for instance \cite[Theorem 9.1.4, vol. 2]{CC00}). Therefore, Lemma \ref{lemma:FP} tells us that every leaf of the induced foliation $\widetilde{\mathcal{F}_i}$ on the universal cover $\widetilde{M_i}$ has compact closure.

Let $N_i$ be as above and consider the cover $q_i : \widehat{N}_i \rightarrow N_i$ defined as $\widetilde{M_i}$ with the solid tori bounded by lifts of $S_i$ drilled, the projection is $q_i = p_i|_{p_i^{-1} (N_i)}$ where $p_i : \widetilde{M}_i \rightarrow M_i$ is the universal cover of $M_i$. We claim that if $p: \widetilde{M} \rightarrow M_1 \# M_2$ is the universal cover of the connected sum, then each connected component of $p^{-1} (N_i)$ is homeomorphic to $\widehat{N}_i$ (notice that each $N_i$ is naturally embedded in $M_1 \# M_2$). 

To see this, it suffice to notice that if $j_i: N_i \subset M_i \hookrightarrow M_1\# M_2$ is the inclusion on the connected sum, then $\ker (j_i)_*$ is isomorphic to $\pi_1 (\widehat{N}_i)$ (by Galois correspondence). By Van-Kampen's theorem $\pi_1 (N_i) \cong \pi_1 (M_i) \# \langle m_i \rangle$, where $m_i$ is the homotopy class of a meridian of the drilled solid torus on $M_i$. The morphism $(j_i)_*$ sends each word $\gamma$ of $\pi_1 (N_i)$ to the word obtained by deleting all appearances of $\langle m_i \rangle$ from $\gamma$, because $i_* : \pi_1 (S) \rightarrow \pi_1 (M_1 \# M_2)$ is zero. Therefore, $\ker (j_i)_*$ is the normalizer of $\langle m_i \rangle$, which is exactly the image of $\pi_1 (\widehat{N}_i) \subset \pi_1 (N_i)$.

To conclude the Lemma, notice that every leaf apart from $S$ is either contained in $N_1$ or $N_2$. In any case, the lift of each of those leaves cannot leave a connected component of $p^{-1} (N_i)$ because its boundary is composed of leaves. The restricted foliation on this connected component is the foliation on $\widehat{N}_i \subset \widetilde{M}_i$. Therefore every leaf has compact closure, and the foliation is uniform.
\end{proof}
\end{lemma}

Lemma \ref{lemma:connectedsum} and Construction \ref{cons:S1S2} gives us the following:

\begin{coro}
\label{coro:connectedsumsS1S2}
Given $k \in \mathbb{N}$, there exists a uniform foliation $\mathcal{F}$ with Reeb components on $\#_{i=1}^k S^1 \times S^2$. Furthermore, this foliation admits a transverse unknot.
\end{coro}

Every foliation on a compact $3-$manifold with a finite fundamental group is uniform because the universal cover is compact. Furthermore, it has Reeb components by Novikov's Theorem. We will see that we can construct such foliations in a way that we fall under the hypothesis of Lemma \ref{lemma:connectedsum}.

For the next construction, we will use the notion of \emph{open book decomposition}. An open book decomposition on a closed $3-$manifold $M$ is a pair $(B,\pi)$, where $B$ is an oriented link on $M$ (called the \emph{binding}) and $\pi : M \setminus B \rightarrow S^1$ is a fibration. We ask the fibers $\pi^{-1} (\cdot)$ to be the interior of compact surfaces whose boundary is $B$; these are called the \emph{pages} of the open book. There is a classical theorem (due to Alexander) saying that every closed $3-$manifold has an open book decomposition (see \cite{Etn06} for a nice introduction to these objects).

\begin{cons}
\label{cons:blocks}
We will construct a (uniform) foliation $\mathcal{F}$ on every closed $3-$manifold $M$ with a finite fundamental group, such that there exists a torus transverse to $\mathcal{F}$ which bounds a solid torus $T$. In particular, every meridian of $\partial T$ is an unknot.  

To see this, we will briefly recall the argument which shows that every closed $3-$manifold admits a foliation. Choose an open book decomposition $(B,\pi)$ on $M$. Let $T$ be a tubular neighborhood of $B$ and identify every connected component of $T$ with $S^1 \times D^2$. With this identification, define a tubular neighborhood $T' \subset T$ whose connected components are $S^1 \times D'^2 \subset S^1 \times D^2$, where $D'^2$ is an open disk properly contained on $D^2$.

On $M \setminus T'$ we have a foliation transverse to the boundary defined by the fibers of the fibration $\pi: M \setminus B \rightarrow S^1$. Spinning this foliation along the boundary and filling the boundary with a Reeb component (see for instance \cite[Example 4.10]{Ca07}), we obtain a foliation on $M$. Notice that $\partial T$ is composed of a torus transversal to the constructed foliation, bounding a solid torus.
\end{cons}

Construction \ref{cons:blocks} and Corollary \ref{coro:connectedsumsS1S2} immediately gives us Theorem \ref{thm:constructions}:

\begin{coro}
Given $k \in \mathbb{N}$ and a choice $M_1,\ldots,M_l$ of $3-$manifolds with finite fundamental group, there exists a uniform foliation with Reeb components on $\left( \#_{i=1}^k S^1 \times S^2 \right) \# \left( \#_{i=1}^l M_i \right)$. 
\end{coro}

We will finish this section by sketching an example of a uniform foliation on the solid torus inducing the trivial foliation by meridians on the boundary. To do so, we will use the following well-known lemma. A visual proof can be found in \cite[Section 5.3]{Har80}:

\begin{lemma}
\label{lemma:borromean}
Identify $\mathbb{T}^3$ as the quotient $\mathbb{R}^3/\mathbb{Z}^3$, and let $T_x, T_y, T_z$ be disjoint tubular neighborhoods of curves on $\mathbb{T}^3$ which are projections of lines parallel to the $x,y$ and $z-$axis respectively. Then if $L$ is the Borromean link on $S^3$ and $N$ is a tubular neighborhood of $L$, there exists a homeomorphism $h: \mathbb{T}^3 \setminus \left( T_x \cup T_y \cup T_z \right) \rightarrow S^3 \setminus N$. Furthermore the homeomorphism takes a meridian of $T_x,T_y$ and $T_z$ to a longitude of a component of the borromean ring, i.e. a curve on a component of $N$ bounding a disk on $S^3$.
\end{lemma}

\begin{figure}[h!]
\begin{center}
\includegraphics[scale=0.8]{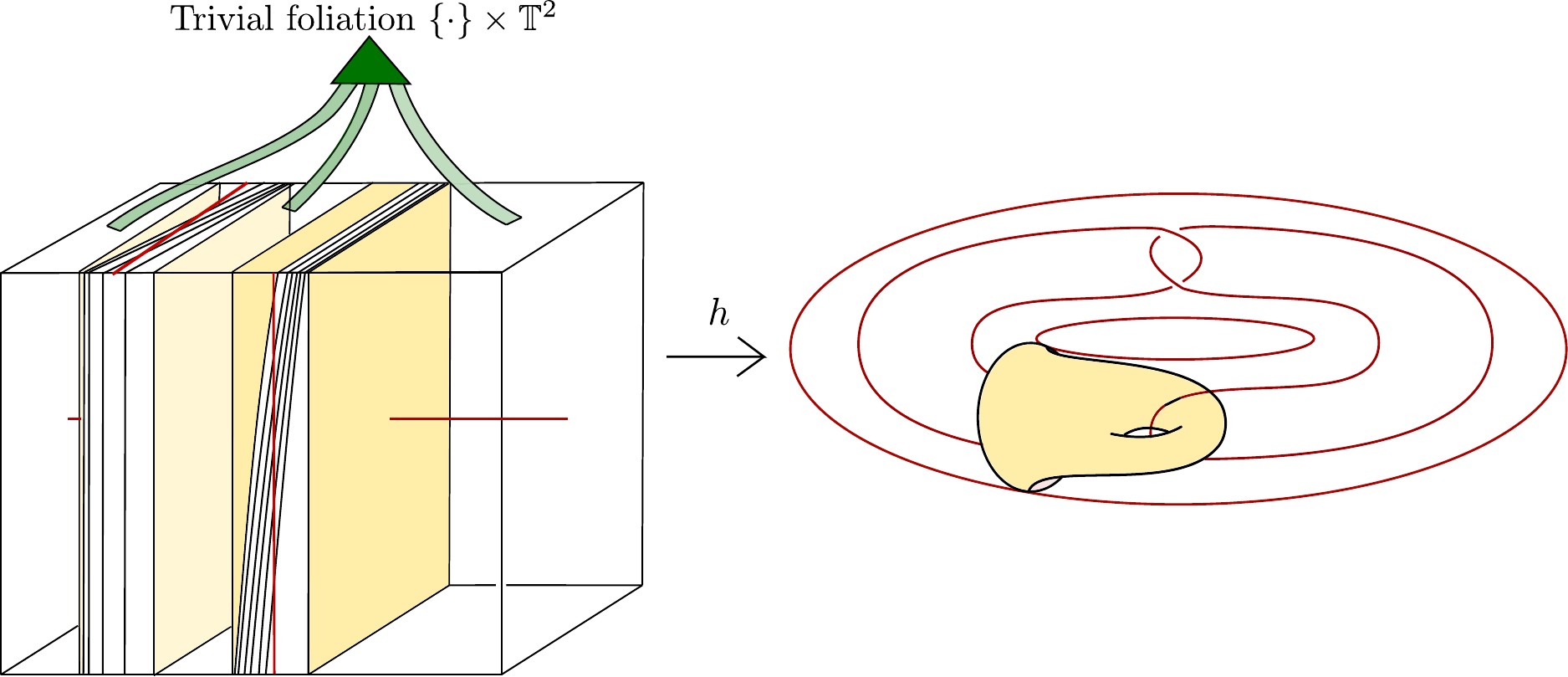}
\caption{At the left: a foliation $\mathcal{F}_1$ on $\mathbb{T}^3$ transverse to curves parallel to the $x,y$ and $z-$axes. At the right: the image of a compact leaf via $h$ on the complement of the Borromean rings.}
\label{fig:borromean}
\end{center}
\end{figure}

\begin{cons}
\label{cons:borromean}
We will construct a uniform foliation with Reeb components on the solid torus $S^1 \times D^2$, inducing a foliation by closed curves of meridional slope on the boundary.

Start with $\mathbb{T}^3$ endowed with a foliation $\mathcal{F}_1$ as the one depicted on Figure \ref{fig:borromean}. This foliation is transversal to three closed curves $c_x,c_y$ and $c_z$ lifting to $\mathbb{R}^3$ as lines parallel to the axis $x,y$ and $z$ respectively. Take tubular neighborhoods $T_x,T_y$ and $T_z$ of $c_x,c_y$ and $c_z$ respectively where the foliation restricts to a trivial foliation by disks.

By Lemma \ref{lemma:borromean}, there is a homeomorphism $h :\mathbb{T}^3 \setminus \left( T_x \cup T_y \cup T_z \right) \rightarrow S^3 \setminus N$, with $N$ a tubular neighborhood of the Borromean link. Notice that $S^3 \setminus N$ is homeomorphic to a solid torus, with the tubular neighborhood of two linked unknots drilled as in the right of Figure \ref{fig:borromean}. To see this, look at the complement of the tubular neighborhood of one link component.

The homeomorphism $h$ sends the foliation on $\mathbb{T}^3$ to a foliation $\mathcal{F}_2$ on the solid torus minus the link. This foliation is transverse to the boundary: it induces on the boundary of the solid torus the trivial foliation by meridians and the trivial foliation by longitudes around the link components. Spiraling around the link in the interior of the solid torus and then adding Reeb components, we obtain a foliation $\mathcal{F}$ on the solid torus, which has the desired behavior on the boundary.

We claim that $\mathcal{F}$ is uniform. To see this take a compact leaf $L$ of $\mathcal{F}_1$ in $\mathbb{T}^3$ intersecting only the curve $c_z$. Without lost of generality, $h(\partial T_z)$ is identified with the boundary of our solid torus with a link drilled. The leaf $L$ does not separate $\mathbb{T}^3 \setminus \left( T_x \cup T_y \cup T_z \right)$. Therefore, its image $h(L)$ is a leaf of the foliation $\mathcal{F}$ which does not separate the solid torus. A non-separating surface on the solid torus must have a nontrivial intersection number with the curves $S^1 \times \{\cdot \}$. This implies that when we lift the foliation to the universal cover, the lifts of this compact leaf bound a fundamental domain for $S^1 \times D^2$. Explicitly, this leaf on the solid torus looks like the yellow surface at the right of Figure \ref{fig:borromean} (see the proof of Lemma \ref{lemma:borromean} in \cite[Section 5.3]{Har80}). We conclude the uniform condition using Lemma \ref{lemma:construction}.
\end{cons}

\section{Proof of Theorem \ref{thm:main}}

In this section, $M$ will be a compact $3-$manifold endowed with a uniform foliation $\mathcal{F}$ with Reeb components. We equip $M$ with a Riemannian metric; this induces on the universal cover $p : \widetilde{M} \rightarrow M$ the pullback metric. Lifting the foliation $\mathcal{F}$ to $\widetilde{M}$ we obtain a foliation $\widetilde{\mathcal{F}}$. By Lemma \ref{lemma:FP}, every leaf of this foliation has compact closure.

We will prove Theorem \ref{thm:main} by contradiction. Suppose that there exists an essential immersion $i: \Sigma \rightarrow M$ from a closed surface $\Sigma$ of genus $g \geq 1$ such that $i(\Sigma)$ does not intersect $\mathcal{R}$, the union of the Reeb components. Up to taking a finite cover, we can assume that $\mathcal{F}$ is oriented and cooriented. Using Theorem \ref{thm:generalposition} we can assume that $i$ is in general position with respect to the foliation. This is because there are immersions $C^0-$close to $i$ which are in general position, therefore if $i(\Sigma)$ is at a finite distance from $\mathcal{R}$, a sufficiently small perturbation of $i$ will also be disjoint from $\mathcal{R}$. The foliation $\mathcal{F}$ induces on $\Sigma$ a singular foliation of Morse-type, which will be denoted as $\mathcal{G}$.

Denote by $q: \mathbb{R}^2 \cong \widetilde{\Sigma} \rightarrow \Sigma$ the universal cover of $\Sigma$. Lifting the immersion we obtain a map $\widetilde{i}: \mathbb{R}^2 \cong \widetilde{\Sigma} \rightarrow \widetilde{M}$. As the immersion is essential, $\widetilde{i}$ is equivariant with respect to the actions of $\pi_1 (\Sigma)$ on $\mathbb{R}^2$ and $\widetilde{M}$. More precisely, for every $\gamma \in \pi_1 (\Sigma)$ we have $\widetilde{i} \circ \gamma = i_* (\gamma) \circ \widetilde{i}$. The immersed plane $\widetilde{i} (\mathbb{R}^2)$ is also in general position with respect to $\widetilde{\mathcal{F}}$. In fact the induced singular foliation on $\mathbb{R}^2$ is exactly the lift of $\mathcal{G}$ to the universal cover, we will denote it by $\widetilde{\mathcal{G}}$. Notice that every leaf of $\widetilde{\mathcal{G}}$ has compact closure because every leaf of $\widetilde{\mathcal{F}}$ has compact closure.

We can classify every closed graph of $\widetilde{\mathcal{G}}$:

\begin{rem}
\label{rem:separatrixesboundsdisk}
Let $S$ be a separatrix of $\widetilde{\mathcal{G}}$ such that its $\alpha$ and $\omega-$limit is a singularity, then $\alpha (S) = \omega (S)$. To see this, let $\{x \}$ and $\{y\}$ be the $\alpha$ and $\omega-$limit respectively. The separatrixes of $\mathcal{G}$ cannot joint two distinct singularities by the second item on the definition of immersions in general position, therefore $y = \gamma x$ for some $\gamma \in \pi_1 (\Sigma)$. If $\gamma$ is not the identity, then $q (S \cup \{x\} \cup \{\gamma x\})$ is a homotopically nontrivial curve. As the immersion is essential, this curve is sent to a homotopically nontrivial curve on $M$ contained on certain leaf of $\mathcal{F}$. This is a contradiction by Lemma \ref{lemma:FP}, therefore $x = y$.
\end{rem}

\begin{figure}[t]
\begin{center}
\includegraphics[scale=1]{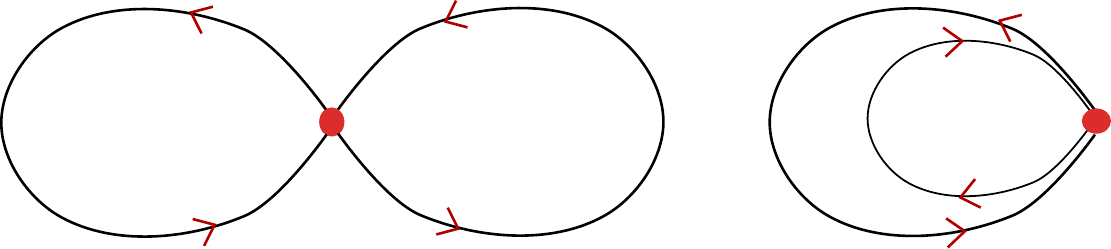}
\caption{The possible configurations for the saddles in $\widetilde{\mathcal{G}}$.}
\label{fig:generalizedclosedleaf}
\end{center}
\end{figure}

We will call a closed leaf or closed graph of $\widetilde{\mathcal{G}}$ which is maximal with respect to inclusion (i.e., is not strictly contained on another closed graph) a \emph{generalized closed leaf}. The last remark tells us that a separatrix limits on at most one singularity. Therefore, if $S$ is a separatrix of a generalized closed leaf limiting on a singularity $x$, then $S \cup \{x\}$ bounds a disk (by the Jordan curve theorem). This leaves us with exactly two possible configurations for a generalized closed leaf containing a singularity. Either the disk bounded by one separatrix contains the disk bounded by the other, or both disks are disjoint; both configurations are depicted in Figure \ref{fig:generalizedclosedleaf}. These disks project homeomorphically to $\Sigma$, because their boundaries project homeomorphically.

Suppose that $D \subset \mathbb{R}^2$ is a disk bounded by a separatrix $S$ of a generalized closed leaf, then there is a well-defined notion of holonomy transport on the side of $D$. More precisely, take a half-open transversal $\tau$ of a point of $S$, with $\tau \subset D$. Start with a point $x$ in $\tau$ and look at the next point of intersection between the leaf passing through $x$ and $\tau$ (with respect to the orientation on $\mathcal{G}$). This defines a map from $\tau$ to itself. Choosing finitely many foliated neighborhoods of $S$ and a model neighborhood of its corresponding singularity, one can see that this map is well defined up to take a smaller half-open transversal. Analogously, we can define the holonomy on the unbounded side of a generalized closed leaf.

\begin{lemma}
\label{lemma:allleavescompact}
Let $\mathcal{F}$ be an oriented and transversely oriented, uniform foliation with Reeb components on a compact $3-$manifold $M$ and let $i: \Sigma \rightarrow M$ be an essential immersion of a surface $\Sigma$ of genus $\geq 1$ in general position with respect to $\mathcal{F}$. Suppose we further assume that $i(\Sigma)$ is disjoint from the set of Reeb components, then every leaf of the induced singular foliation $\mathcal{G}$ on $\Sigma$ is either compact or a separatrix such that its $\alpha$ and $\omega-$limit coincides.
\begin{proof}
First, we prove that no generalized closed leaf of $\widetilde{\mathcal{G}}$ has holonomy in any of its sides. Let $C$ be a generalized closed leaf, if we had nontrivial holonomy on a certain side of $C$, then the curve $q(C) \subset \Sigma$ would be sent via $i$ to a homotopically nontrivial curve on a leaf of $\mathcal{F}$. To see this, if $q(C)$ were sent to a homotopically trivial curve on a leaf of $\mathcal{F}$ then it would have trivial holonomy on $\mathcal{F}$, and therefore trivial holonomy on $\widetilde{\mathcal{G}}$. 

Now we can use standard Haefliger-type arguments: by Remark \ref{rem:separatrixesboundsdisk} $q(C)$ bounds one or two disks which are sent via $i$ to immersed disks on $M$. The curve $i(q(C))$ is non
homotopically trivial on its corresponding leaf of $\mathcal{F}$; therefore, the boundary of at least one of the immersed disks is non homotopically trivial. This implies the existence of a vanishing cycle on at least one of the immersed disks (see for instance \cite[Lemma 9.2.2 and 9.2.4, vol. 2]{CC00}). However, Novikov's Theorem says that if a leaf admits a vanishing cycle, then this leaf is the boundary of a Reeb component (see \cite[Theorem 9.4.1, vol. 2]{CC00}).  Therefore, $i(\Sigma)$ intersects the set of Reeb components, which is a contradiction.

Now we prove that either every leaf of $\mathcal{G}$ is compact or a separatrix limiting on a single singularity. It is enough to see this on $\widetilde{\mathcal{G}}$. Suppose that $L$ is a separatrix of $\widetilde{\mathcal{G}}$ with $\alpha (S) \neq \omega (S)$ or a non-compact leaf. In any case, $\omega (L)$ or $\alpha (L)$ must be a generalized closed leaf or a closed leaf by Remark \ref{rem:separatrixesboundsdisk} and Poincaré-Bendixson theorem. This generalized closed leaf must have nontrivial holonomy in one of its sides, which is a contradiction by the last paragraphs.
\end{proof}
\end{lemma}

It is not hard to see that a singular foliation $\mathcal{G}$ on $\Sigma$ without holonomy must have a generalized closed leaf which is homotopically nontrivial:

\begin{lemma}
\label{lemma:breakingthings}
Let $\Sigma$ be a closed Riemannian surface of genus $g \geq 1$ and $\mathcal{G}$ an oriented singular foliation of Morse-type. Furthermore, suppose that every leaf is compact or a separatrix with the same $\alpha$ and $\omega-$limit. Then some leaf of the lifted singular foliation $\widetilde{\mathcal{G}}$ on $\mathbb{R}^2 \cong \widetilde{\Sigma}$ is unbounded.
\begin{proof}
Assume by contradiction that every leaf of $\widetilde{\mathcal{G}}$ is bounded, therefore has compact closure. Let $\mathcal{C}$ be the set of generalized closed leaves of $\widetilde{\mathcal{G}}$, and define the function $A : \mathcal{C} \rightarrow \mathbb{R}_{\geq 0}$ such that $A(C)$ is the area of the bounded component of the complement of $C \in \mathcal{C}$. Each connected component of the bounded component of an element of $\mathcal{C}$ is a disk that projects homeomorphically to $\Sigma$, therefore its area cannot be greater than the area of $\Sigma$. This shows that the function $A$ is bounded.

We claim that we can achieve the supremum $\sup_{C \in \mathcal{C}} A(C)$. Choose $x_n \in \mathbb{R}^2$ such that there exists a generalized closed leaf $C_n$ with $x_n \in C_n$ and $A(C_n) \to \sup_{C \in \mathcal{C}} A(C)$ an increasing sequence. As the surface $\Sigma$ is compact, up to translate $\{x_n\}$ by $\pi_1 (\Sigma)$ we can assume that the sequence $\{x_n\}_{n\in \mathbb{N}}$ converges to $x \in \mathbb{R}^2$. 

Let $C$ be the generalized closed leaf passing through $x$; we claim that $A(C)$ is the desired supremum. To see this, notice that if $m> n$ then $A(C_n) < A(C_m)$ and therefore, $C_n$ is contained on the bounded component of $C_m$. This implies that for $n$ sufficiently big, $C_n$ is contained on the bounded component of $C$. Otherwise, there would exist some $C_k$ bounding a region containing $C_n$ with $n > k$, which is a contradiction.

Now choose $C$ as in the last paragraph. This leaf cannot have holonomy because every leaf is compact. Therefore we can find a foliated neighborhood of $C$ composed of closed orbits. In particular, some closed leaf bounding $C$ exists, which contradicts that $C$ bounds the largest area.
\end{proof}
\end{lemma}

\begin{figure}[t!]
\begin{center}
\includegraphics[scale=0.6]{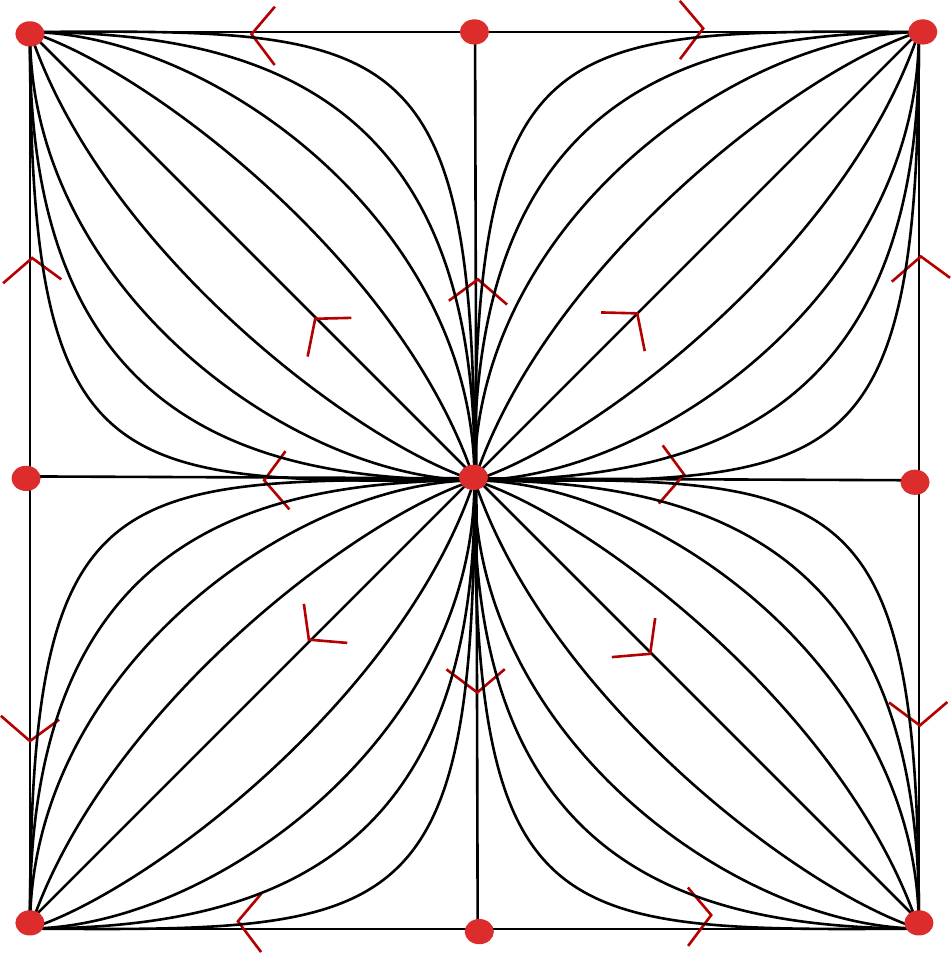}
\caption{Turbulizing around the attracting/repelling fixed points we obtain a uniform foliation on $\mathbb{T}^2$.}
\label{fig:uniformsurface}
\end{center}
\end{figure}

Therefore, Lemma \ref{lemma:allleavescompact} tells us that the induced foliation on the surface $\Sigma$ satisfies the hypothesis of Lemma \ref{lemma:breakingthings}. However, this lemma contradicts the uniformity of the foliation $\mathcal{F}$. This contradiction concludes the proof of Theorem \ref{thm:main}.

\begin{rem}
Suppose that $\mathcal{F}$ is an oriented and transversely oriented, uniform foliation with Reeb components on a compact $3-$manifold $M$, and let $i: \Sigma \rightarrow M$ be an essential immersion in general position as before. We just proved that $i(\Sigma)$ must intersect the set of Reeb components. Thus, the induced singular foliation on $\Sigma$ is uniform and necessarily contains vanishing cycles. 

There exist singular foliations on surfaces satisfying these conditions, see for example Figure \ref{fig:uniformsurface}.\footnote{We thank Sergio Fenley for pointing out this example.} It is not obvious how the singular foliation induced on $\Sigma$ would help in order to answer Question \ref{ques}. We believe that an extra argument is needed.
\end{rem}

\bibliographystyle{alpha}

\end{document}